\documentclass[twoside]{article} 
\date{} 
\title{The asymptotic expansion  of a Mathieu-exponential series}
\author{\sc R. B.\ Paris \\
{\em Division of Computing and Mathematics,} \\
{\em Abertay University, Dundee DD1 1HG, UK}}
\begin{document}
\def\f#1#2{\mbox{${\textstyle \frac{#1}{#2}}$}}
\def\dfrac#1#2{\displaystyle{\frac{#1}{#2}}}
\def\boldal{\mbox{\boldmath $\alpha$}}
\newcommand{\bee}{\begin{equation}}
\newcommand{\ee}{\end{equation}}
\newcommand{\lam}{\lambda}
\newcommand{\ka}{\kappa}
\newcommand{\al}{\alpha}
\newcommand{\ba}{\beta}
\newcommand{\la}{\lambda}
\newcommand{\ga}{\gamma}
\newcommand{\eps}{\epsilon}
\newcommand{\fr}{\frac{1}{2}}
\newcommand{\fs}{\f{1}{2}}
\newcommand{\g}{\Gamma}
\newcommand{\br}{\biggr}
\newcommand{\bl}{\biggl}
\newcommand{\ra}{\rightarrow}
\newcommand{\gtwid}{\raisebox{-.8ex}{\mbox{$\stackrel{\textstyle >}{\sim}$}}}
\newcommand{\ltwid}{\raisebox{-.8ex}{\mbox{$\stackrel{\textstyle <}{\sim}$}}}
\renewcommand{\topfraction}{0.9}
\renewcommand{\bottomfraction}{0.9}
\renewcommand{\textfraction}{0.05}
\newcommand{\mcol}{\multicolumn}
\date{}
\maketitle
\pagestyle{myheadings}
\markboth{\hfill \sc R. B.\ Paris  \hfill}
{\hfill \sc A Mathieu-exponential series\hfill}
\begin{abstract}
We consider the asymptotic expansion of the functional series 
\[S_{\mu}^\pm(a;\lambda)=\sum_{n=0}^\infty \frac{(\pm 1)^n e^{-\lambda n}}{(n^2+a^2)^\mu}\]
for $\lambda>0$ and $\mu\geq0$ as $|a|\to \infty$ in the sector $|\arg\,a|<\pi/2$. The approach employed
consists of expressing $S_{\mu}^\pm(a;\lambda)$ as a contour integral combined with suitable deformation of the integration path.
Numerical examples are provided to illustrate the accuracy of the various expansions obtained.

\vspace{0.3cm}

\vspace{0.1cm}
 
\end{abstract}

\vspace{0.3cm}

\noindent $\,$\hrulefill $\,$

\vspace{0.3cm}

\begin{center}
{\bf 1.\ Introduction}
\end{center}
\setcounter{section}{1}
\setcounter{equation}{0}
\renewcommand{\theequation}{\arabic{section}.\arabic{equation}}
The functional series
\bee\label{e11}
\sum_{n=1}^\infty\frac{n}{(n^2+a^2)^\mu}
\ee
is known as a Mathieu series \cite{M}, which originally arose (in the case $\mu=2$) in problems dealing with the elasticity of solid bodies.  
The asymptotic expansion for large $a$ of more general functional series of this type has been discussed in \cite{P} and \cite{Z}. 
More recently, Gerhold and Tomovski \cite{G} extended the asymptotic study of (\ref{e11}) by introducing the factor $z^n$, where $|z|\leq 1$. From this result they were able to deduce, in particular, the large-$a$ expansions of the trigonometric Mathieu series
\[\sum_{n=1}^\infty \frac{n \sin nx}{(n^2+a^2)^\mu},\qquad\sum_{n=1}^\infty \frac{n \cos nx}{(n^2+a^2)^\mu}.\]

Subsequently, the above trigonometric series were generalised in the form \cite{P3}
\bee\label{e12}
\sum_{n=1}^\infty \frac{n^\gamma {\cal C}_\nu(bn/a)}{(n^2+a^2)^\mu}\qquad (b>0),
\ee
where ${\cal C}_\nu$ denotes the oscillatory Bessel functions $J_\nu(x)$ and $Y_\nu(x)$ with argument proportional to $n/a$, and their large-$a$ asymptotics determined. 
In addition, this last study also considered the inclusion of the modified Bessel function $K_\nu(x)$ of similar argument,
which contains the decaying exponential as a special case when $\nu=\fs$, since $$K_{1/2}(x)=(\pi/2x)^{1/2} e^{-x}.$$ In \cite{P4}, the case where the additional term is a Gaussian exponential, namely the series
\[\sum_{n=1}^\infty \frac{n^\gamma e^{-\la n^2/a^2}}{(n^2+a^2)^\mu}\qquad (\la>0),\]
has been considered. For even integer values of the parameter $\gamma$ it is found that the large-$a$ asymptotic expansion consists of an algebraic expansion with a finite number of terms together with a sequence of increasingly subdominant exponentially small contributions. This situation is analogous to the well-known Poisson-Jacobi transformation (corresponding to $\mu=\gamma=0$) given by \cite[p.~124]{WW}
\[\sum_{n=1}^\infty e^{-\la n^2/a^2}=\frac{a}{2}\sqrt{\frac{\pi}{\la}}-\frac{1}{2}+a\sqrt{\frac{\pi}{\la}} \sum_{n=1}^\infty e^{-\pi^2n^2a^2/\la}.\]

The asymptotic expansion we consider here is the alternating Mathieu series coupled with a decaying exponential (depending linearly on the summation index $n$) of the form
\bee\label{e13}
S_{\mu}^\pm(a;\lambda):=\sum_{n=0}^\infty\frac{(\pm 1)^{n} e^{-\la n}}{(n^2+a^2)^\mu}\qquad (\mu\geq0,\ \lambda>0)
\ee
for $|a|\to\infty$ in the sector $|\arg\,a|<\fs\pi$. In \cite[\S 5]{P3}, the related series with the exponential factor $e^{-bn/a}$, resulting from (\ref{e12}) with ${\cal C}_\nu(nb/a)=K_\nu(bn/a)$ when $\nu=\gamma=\fs$, was shown to possess the asymptotic expansion
\[\sum_{n=1}^\infty\frac{(\pm 1)^n e^{-b n/a}}{(n^2+a^2)^\mu}-a^{2\mu-1} \bl(\frac{\pi}{2b}\br)^{\!\!1/2}\bl\{\!\!\!\begin{array}{c} J_\mu(a; b/a)\\0\end{array}\!\!\!\br\}\sim\sqrt{\frac{2b}{\pi a}}\,\{R^\pm(a;\fs)+R^\pm(a;-\fs)\}\]
as $|a|\to\infty$ in $|\arg\,a|<\fs\pi$, where the quantity $J_\mu(a;b/a)$ is given in (\ref{e31}) and
\[R^\pm(a;w):=\frac{a^{-w-2\mu}b^w}{2^{w+1}} \g(-w) \sum_{k=0}^\infty\frac{(-1)^k (\mu)_k}{k! a^{2k}} Z^\pm(-\omega_k) F_k(w,b).\]
Here $Z^+(x):=\zeta(x)$, $Z^-(x):=(1-2^{1-x}) \zeta(x)$, $\omega_k:=\fs+w+2k$, with $\zeta(s)$ being the Riemann zeta function, and $F_k(w,b)$ are polynomials expressed as a terminating hypergeometric series defined by
\[F_k(w,b):={}_1F_2(-k;\fs+w,1-\mu-k; -b^2/4).\]
However, if we set $b=\la a$, where $\la>0$ is finite, to obtain a series equivalent to that in (\ref{e13}), it is seen that the polynomials $F_k(w,b)=F_k(w,\la a)$ with the consequence that the formal series $R(a;\pm \fs)$ lose
their asymptotic character. 

Rather than adopting a Mellin transform approach used in \cite{P3}, we express $S_\mu^\pm(a;\la)$ as a contour integral combined with suitable integration path deformation. Such an approach has been employed by Olver in his well-known book \cite[p.~303]{Olv}, who showed that in the particular case $\la=0$ 
\bee\label{e14}
\sum_{n=0}^\infty \frac{(-1)^n}{(n^2+a^2)^\mu}=\frac{1}{2a^{2\mu}}+\frac{2^{\frac{3}{2}-\mu}\sqrt{\pi}}{a^{2\mu-1} \g(\mu)} \sum_{k=0}^\infty \frac{K_{\frac{1}{2}-\mu}((2k+1)\pi a)}{((2k+1)\pi a)^{\frac{1}{2}-\mu}}
\ee
for $\mu>0$.
For large values of $a$ in $|\arg\,a|<\fs\pi$, the sum of Bessel functions decays exponentially fast.
Series of the form (\ref{e14}) have arisen in aerodynamic interference calculations; see \cite{Olv49}.
We shall find that when $\la>0$ there is an analogous sum of $K$-Bessel functions (with complex argument) together with an additional term possessing an algebraic asymptotic expansion when $a$ is large.

We note at this point the special evaluations for non-negative integer $\mu$ given by
\bee\label{e15}
S_0^\pm(a;\la)=\frac{e^\la}{e^\la\mp 1}
\ee
and 
\begin{eqnarray}
S_1^\pm(a;\la)\!\!&=&\!\!\frac{{\cal F}_1}{2a^2},\qquad S_2^\pm(a;\la)=\frac{1}{4a^4}({\cal F}_1+{\cal F}_2),\nonumber\\
S_3^\pm(a;\la)\!\!&=&\!\!\frac{1}{16a^6} (3{\cal F}_1+3{\cal F}_2+2{\cal F}_3),\nonumber\\
S_4^\pm(a;\la)\!\!&=&\!\!\frac{1}{32a^8}(5{\cal F}_1+5{\cal F}_2+4{\cal F}_3+2{\cal F}_4),\nonumber\\
S_5^\pm(a;\la)\!\!&=&\!\!\frac{1}{256a^{10}}(35{\cal F}_1+35{\cal F}_2+30{\cal F}_3+20{\cal F}_4+8{\cal F}_5),
 \dots\, ,\label{e16}
\end{eqnarray}
where
\[{\cal F}_n=F_n(a)+F_n(-a),\quad F_n(a)={}_{n+1}F_n\bl(\!\!\begin{array}{c}1, ia, \ldots , ia\\1+ia, \ldots , 1+ia\end{array};\,\pm e^{-\la}\br)\]
and ${}_{n+1}F_n$ is the generalised hypergeometric function.

\vspace{0.6cm}

\begin{center}
{\bf 2.\ Derivation of the expansion for $S_\mu^-(a;\la)$}
\end{center}
\setcounter{section}{2}
\setcounter{equation}{0}
\renewcommand{\theequation}{\arabic{section}.\arabic{equation}}
We first consider the sum $S_\mu^-(a;\la)$.
Following \cite[p.~303]{Olv}, we have
\bee\label{e21}
S_\mu^-(a;\la)=\frac{1}{2i}\int_\infty^{(0+)} \frac{e^{-\la t}}{\sin \pi t}\,\frac{dt}{(a^2+t^2)^\mu},\qquad (|\arg\,a|<\fs\pi),
\ee
where the integration path encloses only the poles of the integrand situated at $t=0, 1, 2, \ldots\ $. Let us first consider the case $a>0$. The integrand then has branch points at $t=\pm ia$ with cuts along the imaginary axis
emanating from these points and passing to infinity. Since $\la>0$, the loop contour may be deformed to coincide with the imaginary axis between $[-i(a-\rho), i(a-\rho)]$ (with an indentation of radius $\rho$ around the origin on the left-hand side) together with the portions of the imaginary axis situated above $ia$ and below $-ia$
on the right-hand side of the cuts; see Fig.~1. Provided $\mu<1$, the contribution from the indentations of radius $\rho$ round the branch points vanishes as $\rho\to 0$.

\begin{figure}[t]
\centering
\begin{picture}(200,200)(0,0)
\put(0,100){\line(1,0){150}}
\put(20,155){\line(0,1){40}}
\put(18,155){\line(0,1){40}}
\put(17.5,155){\line(0,1){40}}
\put(20,105){\line(0,1){40}}
\put(20,100){\oval(10,10)[l]}
\put(20,150){\oval(10,10)[r]}
\put(20,120){\line(1,1){10}}
\put(20,120){\line(-1,1){10}}
\put(20,95){\line(0,-1){40}}
\put(20,45){\line(0,-1){40}}
\put(18,45){\line(0,-1){40}}
\put(17.5,45){\line(0,-1){40}}
\put(20,50){\oval(10,10)[r]}
\put(20,100){\circle*{3}}
\put(20,150){\circle*{3}}
\put(20,50){\circle*{3}}
\put(50,100){\circle*{3}}
\put(80,100){\circle*{3}}
\put(110,100){\circle*{3}}
\put(140,100){\circle*{3}}
\put(25,88){0}
\put(50,88){1}
\put(80,88){2}
\put(110,88){3}
\put(140,88){4}
%
\put(5,148){$ia$}
\put(0,48){$-ia$}
\end{picture}
\caption{\small{The deformed integration path in the $t$-plane. The heavy lines emanating from the points $t=\pm ia$ denote branch cuts. The arrow denotes the direction of integration.}}
\end{figure}
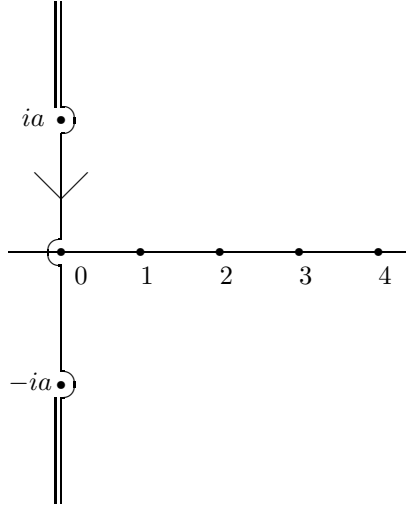

\vspace{0.2cm}
\newpage
\noindent{\bf 2.1\ \ The contribution between the branch points}
\vspace{0.2cm}

\noindent
We first deal with the integral between the branch points which we denote by $I_1$. We have
\[I_1=-\frac{a^{1-2\mu}}{2i} \int_{-1+\rho}^{1-\rho}\frac{e^{-\la iat}}{\sinh \pi at}\,\frac{dt}{(1-t^2)^\mu},\]
where the integration path is indented to the left at the origin, so that
\bee\label{e20}
I_1=\frac{1}{2a^{2\mu}}+H_\mu^-(a;\la), \qquad H_\mu^-(a;\la):=a^{1-2\mu}\int_0^1\frac{\sin \la at}{\sinh \pi at}\,\frac{dt}{(1-t^2)^\mu}
\ee
when $\rho\to 0$. Now
\[\frac{\sin \la x}{\sinh \pi x}=\frac{\la}{\pi}\sum_{k=0}^\infty (-1)^k A_k x^{2k}\qquad (|x|<\infty),\]
where
\begin{eqnarray*}
A_0&=&1,\quad A_1=\frac{1}{6} (\la^2+\pi^2),\quad A_2=\frac{1}{360} (3\la^4+10\la^2\pi^2+7\pi^4),\\
A_3&=&\frac{1}{15120} (3\la^6+21\la^4\pi^2+49\la^2\pi^4-31\pi^6),\\
A_4&=&\frac{1}{1814400} (5\la^8+60\la^6\pi^2+294\la^4\pi^4+620\la^2\pi^6+381\pi^8), \ldots\ .
\end{eqnarray*}
Then we find
\[H_\mu^-(a;\la)=\frac{\la a^{1-2\mu}}{\pi} \sum_{k=0}^\infty (-1)^k A_k a^{2k}\int_0^1 \frac{t^{2k}}{(1-t^2)^\mu}\,dt\]
\bee\label{e22}
=\frac{\la a^{1-2\mu}}{2\pi} \g(1-\mu) \sum_{k=0}^\infty \frac{(-1)^k A_k \g(k+\fs)}{\g(k+\frac{3}{2}-\mu)}\,a^{2k} \qquad (\mu<1).
\ee
This form of expansion is suitable for computation when $a$ is not large.

To deal with the situation when $a\to\infty$, we expand the factor $(1-t^2)^{-\mu}$ appearing in (\ref{e20}) by the binomial theorem to obtain
\[H_\mu^-(a;\la)=a^{1-2\mu} \sum_{k=0}^\infty \frac{(\mu)_k}{k!} \int_0^1 t^{2k} \frac{\sin \la at}{\sinh \pi at}\,dt.\]
Since $\sinh \pi at$ decays exponentially for $t$ bounded away from the origin, we can replace the upper limit 1 by $\infty$, thereby introducing an error of O($e^{-\pi a})$. From the result \cite[(3.524.1)]{GR}, we have
\[B_k:=\int_0^\infty \tau^{2k}\,\frac{\sin \la\tau}{\sinh \pi\tau}\,d\tau=\frac{i(2k)!}{(2\pi )^{2k+1}}\bl\{\zeta\bl(2k+1,\frac{1}{2}+\frac{i\la}{2\pi}\br)-\zeta\bl(2k+1,\frac{1}{2}-\frac{i\la}{2\pi}\br)\br\}\]
where $\zeta(a,s)=\sum_{n=0}^\infty (n+a)^{-s}$ is the Hurwitz zeta function. Using the fact that $\zeta(2k+1, z)=-\psi^{(2k)}(z)/(2k)!$, where $\psi(z)$ is the logarithmic derivative of the gamma function, and \cite[(5.15.6)]{DLMF}
\bee\label{epsi}
\psi^{(2k)}(1-z)-\psi^{(2k)}(z)=\pi \bl(\frac{d}{dz}\br)^{\!\!2k}\!\! \cot \pi z,
\ee
we find that
\[B_k=\frac{i}{(2\pi)^{2k+1}}\bl\{\psi^{(2k)}\bl(\frac{1}{2}- \frac{i\la}{2\pi}\br)-\psi^{(2k)}\bl(\frac{1}{2}+\frac{i\la}{2\pi}\br)\br\}=\frac{(-1)^k}{2^{2k+1}} \bl(\frac{d}{dx}\br)^{\!\!2k} \tanh x\br|_{x=\la/2}.\]
The first few coefficients $B_k$ are given by:
\begin{eqnarray*}
B_0&=&\frac{1}{2} \tanh x,\quad B_1=\frac{\sinh x}{4\cosh^3 x},\quad B_2=\frac{\sinh x}{4\cosh^5 x}(2-\sinh^2 x),\\
B_3&=&\frac{\sinh x}{8\cosh^7x}\,(17-26\sinh^2x+2\sinh^4x),\\
B_4&=&\frac{\sinh x}{4\cosh^9x}\,(62-192\sinh^2x+60\sinh^4x-\sinh^6x), \ldots\ ,
\end{eqnarray*}
where $x=\la/2$.
Hence we obtain the expansion
\bee\label{e24}
H_\mu^-(a;\la)\sim a^{-2\mu} \sum_{k=0}^\infty \frac{(\mu)_k B_k}{k!\, a^{2k}}
\ee
as $a\to+\infty$.
\vspace{0.2cm}

\noindent{\bf 2.2\ \ The contribution beyond the branch points}
\vspace{0.2cm}

\noindent
Consider the contribution to $S_\mu^-(a;\la)$ resulting from the part of the integration path situated in the interval $[ia, i\infty)$. We have
\[I_2=-\frac{1}{2i}\int_{ia}^{\infty i}\frac{e^{-\la t}}{\sin \pi t}\,\frac{dt}{(a^2+t^2)^\mu}=
-\frac{a^{1-2\mu}e^{-\pi i\mu}}{2i} \int_1^\infty \frac{e^{-\la iat}}{\sinh \pi at}\,\frac{dt}{(t^2-1)^\mu}.\]
Now on $t\in[1,\infty)$
\[\mbox{cosech}\, \pi at=2e^{-\pi at} (1-e^{-2\pi at})^{-1}=2\sum_{k=0}^\infty e^{-(2k+1)\pi at},\]
so that substitution in $I_2$ followed by term-by-term integration yields
\bee\label{e24a}
I_2=ie^{-\pi i\mu} a^{1-2\mu} \sum_{k=0}^\infty \int_1^\infty \frac{e^{-X_k t}}{(t^2-1)^\mu}\,dt,\qquad X_k:=(2k+1)\pi a+\la ia
\ee
\[=\frac{i\sqrt{\pi}\,e^{-\pi i\mu}}{\g(\mu) \sin \pi\mu}\,(2a^2)^{\frac{1}{2}-\mu} \sum_{k=0}^\infty \frac{K_{\frac{1}{2}-\mu}(X_k)}{(X_k)^{\frac{1}{2}-\mu}},\]
where the integral has been evaluated in terms of a modified Bessel function by \cite[p.~254]{Olv} when $\mu<1$.

The contribution from the path $[-ia, -i\infty)$ is the conjugate of the above expression. Hence we find
\bee\label{e25}
I_2+{\overline I}_2=T_\mu(a;\la):=\frac{2^{\frac{3}{2}-\mu}\sqrt{\pi}}{a^{2\mu-1} \g(\mu)} \sum_{k=0}^\infty \frac{\sin (\pi \mu-\theta_k)}{\sin \pi\mu} \bl|\frac{K_{\frac{1}{2}-\mu}(X_k)}{(X_k)^{\frac{1}{2}-\mu}}\br|,
\ee
where
\[\theta_k:=\arg \bl(X_k^{\mu-\frac{1}{2}}K_{\frac{1}{2}-\mu}(X_k)\br).\]

Collecting together the results in (\ref{e20}) and (\ref{e25}), we finally obtain the expansion
\bee\label{e26}
S_\mu^-(a;\la)=\frac{1}{2a^{2\mu}}+H_\mu^-(a;\la)+\frac{2^{\frac{3}{2}-\mu}\sqrt{\pi}}{a^{2\mu-1} \g(\mu)} \sum_{k=0}^\infty \frac{\sin (\pi \mu-\theta_k)}{\sin \pi\mu} \bl|\frac{K_{\frac{1}{2}-\mu}(X_k)}{(X_k)^{\frac{1}{2}-\mu}}\br|,
\ee
where $X_k$ is defined in (\ref{e24a}). 
When $\la=0$, we have $H_\mu^-(a;\la)\equiv 0$, $\theta_k\equiv 0$ and (\ref{e26}) reduces to Olver's expansion stated in (\ref{e14}). However, the result (\ref{e26}) has been established provided $\mu<1$. It does not appear possible to extend the validity of (\ref{e26}) to $\mu\geq 1$ by analytic continuation, as was the case in Olver's treatment when $\la=0$. When $a$ is large, $S_\mu^-(a;\la)$ is seen to consist of an algebraic expansion given in (\ref{e24}) together with an exponentially small contribution from the infinite sum of modified Bessel functions.

When $a$ is allowed to take on complex values in the sector $|\arg\,a|<\fs\pi$, the branch points $\pm ia$ move off the imaginary axis to the points $ae^{(\phi\pm\frac{1}{2}\pi)i}$, where $\phi=\arg\,a$, as indicated in \cite[p.~303]{Olv}.  The analysis in this case follows that given above and we conclude that the expansion (\ref{e26}) holds for $|a|\to\infty$ in the sector $|\arg\,a|<\fs\pi$. 

\vspace{0.6cm}

\begin{center}
{\bf 3.\ Derivation of the expansion for $S_\mu^+(a;\la)$}
\end{center}
\setcounter{section}{3}
\setcounter{equation}{0}
\renewcommand{\theequation}{\arabic{section}.\arabic{equation}}
Since the residue of $\cot \pi t$ at any integer $n$ is $1/\pi$, we have
\[S_\mu^+(a;\la)=\sum_{n=0}^\infty \frac{e^{-\la n}}{(n^2+a^2)^\mu}=\frac{1}{2i}\int_\infty^{(0+)} \frac{\cot \pi t\,e^{-\la t}}{(t^2+a^2)^\mu}\,dt\]
where the path encloses only the poles at $t=0, 1, 2, \ldots\,$. The integration path can be deformed as illustrated in Fig.~1; the upper and lower parts of this path are denoted by $C_1$ and $C_2$, respectively.
We have (see \cite[p.~305]{Olv})
\[\frac{\cot \pi t}{2i} =\left\{\begin{array}{ll}\!\!\!\!-\dfrac{1}{2}+\dfrac{e^{\pi it}}{2i \sin \pi t} & (\mbox{on}\ \ C_1)\\
\\
\dfrac{1}{2}+\dfrac{e^{-\pi it}}{2i \sin \pi t} & (\mbox{on}\ \ C_2),\end{array}\right.\]
so that
\[S_\mu^+(a;\la)=J_\mu(a;\la)-\frac{1}{2i}\int_{C_1}\frac{e^{-\la t+\pi it}}{\sin \pi t\,(t^2+a^2)^\mu}\,dt+\frac{1}{2i}\int_{C_2}\frac{e^{-\la t-\pi it}}{\sin \pi t\,(t^2+a^2)^\mu}\,dt.
\]
The integral $J_\mu(a;\la)$ is given by
\bee\label{e31}
J_\mu(a;\la):=
\int_0^\infty\frac{e^{-\la t}}{(t^2+a^2)^\mu}\,dt=\frac{\sqrt{\pi} a^{1-2\mu} \g(1-\mu)}{2(\fs\la a)^{\frac{1}{2}-\mu}}\,{\bf K}_{\frac{1}{2}-\mu}(\la a),
\ee
where ${\bf K}_\nu(x):={\bf H}_\nu(x)-Y_\nu(x)$ is the Struve function defined in \cite[(11.5.2)]{DLMF}.

The contribution between the branch points $\pm ia$ is
\[\frac{1}{2a^{2\mu}}+H^+_\mu(a;\la),\]
where
\[H^+_\mu(a;\la)=a^{1-2\mu}\int_0^1 e^{-\pi au}\,\frac{\sin \la au}{\sinh \pi au}\,\frac{du}{(1-u^2)^\mu}.\]
This integral may be evaluated as in Section 2 by expanding $\sin \la au/\sinh \pi au$. However the resulting integrals are expressible in terms of two ${}_2F_2$ hypergeometric functions and so will not be presented here. Our main interest is the estimation of $H^+_\mu(a;\la)$ for large $a$. Proceeding as in Section 2 we find that
\[H^+_\mu(a;\la)\sim a^{1-2\mu} \sum_{k=0}^\infty\frac{(\mu)_k}{k!} \int_0^\infty e^{-\pi au}\,\frac{\sin \la au}{\sinh \pi au}\,u^{2k}\,du\]
\[
=a^{1-2\mu} \sum_{k=0}^\infty\frac{(\mu)_k}{k! (\pi a)^{2k+1}}\int_0^\infty e^{-w}\,\frac{\sin \la w/\pi}{\sinh w}\,w^{2k}dw.\]

For $k\geq 1$, we have
\[\frac{1}{2i} \int_0^\infty \frac{e^{-w(1\pm i\la/\pi)}}{\sinh w}\,w^{2k}dw=\frac{-i (2k)!}{2^{2k+1}}\,\zeta\bl(2k+1,1\pm \frac{i\la}{2\pi}\br)\]
upon use of the result \cite[(3.552.1)]{GR}
\[\int_0^\infty \frac{x^{\al-1}e^{-\beta x}}{\sinh x}\,dx=2^{1-\al} \g(\al)\,\zeta(2k+1,\fs+\fs\beta)\]
valid for $\Re (\al)>1$, $\Re (\beta)>-1$. Then, we define the coefficients ${\hat B}_k$ by
\[{\hat B}_k:=\frac{i(-1)^{k-1}(2k)!}{(2\pi)^{2k+1}}\bl\{\zeta\bl(2k+1,1-\frac{i\la}{2\pi}\br)-\zeta\bl(2k+1,1+\frac{i\la}{2\pi}\br)\br\}\]
\[=\frac{i(-1)^{k-1}}{(2\pi)^{2k+1}}\bl\{\psi^{(2k)}\bl(1+\frac{i\la}{2\pi}\br)-\psi^{(2k)}\bl(-\frac{i\la}{2\pi}\br)-\frac{i(-1)^k (2k)!}{(\la/2\pi)^{2k+1}}\br\},\]
where we have employed the result  $\psi^{(2k)}(1+z)=\psi^{(2k)}(z)+(2k)!/z^{2k+1}$ \cite[(15.5.5)]{DLMF}.
Application of (\ref{epsi}) then yields
\[{\hat B}_k=2^{-2k-1}\bl\{\bl(\frac{d}{dx}\br)^{\!\!2k} \coth x-\frac{(2k)!}{x^{2k+1}}\br\}_{x=\la/2}.\]
The first few coefficients are given by:
\[{\hat B}_0=\frac{1}{2}\bl(\coth x-\frac{1}{x}\br),\quad {\hat B}_1=\frac{1}{4}\bl(\frac{\cosh x}{\sinh^3 x}-\frac{1}{x^3}\br),\quad {\hat B}_2=\frac{1}{4}\bl(\frac{\cosh x}{\sinh^5x}\,(2+\cosh^2x)-\frac{3}{x^5}\br),\]
\[{\hat B}_3=\frac{1}{8}\bl(\frac{\cosh x}{\sinh^7x}\,(17+26\cosh^2x+2\cosh^4x)-\frac{45}{x^7}\br),\]
\[{\hat B}_4=\frac{1}{4}\bl(\frac{\cosh x}{\sinh^9x}\,(62+192\cosh^2x+60\cosh^4x+\cosh^6x)-\frac{315}{x^9}\br), \ldots\ ,\]
where $x=\la/2$.
Then  we have the expansion
\bee\label{e32}
H^+_\mu(a;\la)\sim a^{-2\mu} \sum_{k=0}^\infty \frac{(-1)^k(\mu)_k{\hat B}_k}{k! a^{2k}}
\ee
as $a\to\infty$.

Finally, the contribution from $[ia, \infty i)$ is
\[I_2=-e^{-\pi i\mu} \frac{a^{1-2\mu}}{2i} \int_1^\infty \frac{e^{-(\pi+i\la)at}}{\sinh \pi at}\,\frac{dt}{(t^2-1)^\mu}=ie^{-\pi i\mu} a^{1-2\mu} \sum_{k=0}^\infty \int_1^\infty \frac{e^{-{\hat X}_kt}}{(t^2-1)^\mu}\,dt\]
\[=\frac{i\sqrt{\pi} e^{-\pi i\mu}}{(2a^2)^{\mu-\frac{1}{2}} \g(\mu) \sin \pi\mu}\,\sum_{k=0}^\infty \frac{K_{\frac{1}{2}-\mu}({\hat X}_k)}{{\hat X}_k^{\frac{1}{2}-\mu}},\]
where 
\bee\label{e33a}
{\hat X}_k:=(2k+2)\pi a+\la ia.
\ee
The contribution from the path $[-ia, -i\infty)$ is the conjugate of the above expression. Hence we find
\bee\label{e33}
I_2+{\overline I}_2={\hat T}_\mu(a;\la):=\frac{2^{\frac{3}{2}-\mu}\sqrt{\pi}}{a^{2\mu-1} \g(\mu)} \sum_{k=0}^\infty \frac{\sin (\pi \mu-{\hat \theta}_k)}{\sin \pi\mu} \bl|\frac{K_{\frac{1}{2}-\mu}({\hat X}_k)}{({\hat X}_k)^{\frac{1}{2}-\mu}}\br|,
\ee
where
\[{\hat \theta}_k:=\arg \bl({\hat X}_k^{\mu-\frac{1}{2}}K_{\frac{1}{2}-\mu}({\hat X}_k)\br).\]

Collecting together the results in (\ref{e31}) -- (\ref{e33}), we finally obtain the expansion
\bee\label{e34}
S_\mu^+(a;\la)=\frac{1}{2a^{2\mu}}+J_\mu(a;\la)+H_\mu^+(a;\la)+\frac{2^{\frac{3}{2}-\mu}\sqrt{\pi}}{a^{2\mu-1} \g(\mu)} \sum_{k=0}^\infty \frac{\sin (\pi \mu-{\hat \theta}_k)}{\sin \pi\mu} \bl|\frac{K_{\frac{1}{2}-\mu}({\hat X}_k)}{({\hat X}_k)^{\frac{1}{2}-\mu}}\br|
\ee
valid when $\mu<1$, where ${\hat X}_k$ is defined in (\ref{e33a}). The exponentially small sum of Bessel functions is seen to be $O(a^{1-2\mu} e^{-2\pi a})$ for large $a$, which is smaller than the equivalent sum arising in $S_\mu^-(a;\la)$ of $O(a^{1-2\mu} e^{-\pi a})$.

When $a$ is complex satisfying $|\arg\,a|<\fs\pi$, the analysis follows the same procedure and we conclude that the expansion (\ref{e34}) holds for $|a|\to\infty$ in the sector $|\arg\,a|<\fs\pi$.
\vspace{0.6cm}

\begin{center}
{\bf 4.\ Numerical verification}
\end{center}
\setcounter{section}{4}
\setcounter{equation}{0}
\renewcommand{\theequation}{\arabic{section}.\arabic{equation}}
To demonstrate the validity of the sum $T_\mu(a;\la)$ involving the modified Bessel functions appearing in (\ref{e26}), we define
\[{\cal S}:=S_\mu^-(a;\la)-H_\mu^-(a;\la)-\frac{1}{2a^{2\mu}},\]
where the contribution $H_\mu^-(a;\la)$ is either computed for small $a$ by (\ref{e22}), or for larger $a$ from (\ref{e20}) using high-precision numerical evaluation of the integral. The value of ${\cal S}$ so obtained is then compared to $T_\mu(a;\la)$ defined in (\ref{e25}). For example, when $a=3$, $\la=1$ and $\mu=\fs$, we find
${\cal S}\doteq-6.35783\,82469\,54\times 10^5$ with $T_\mu(a;\la)\doteq-6.35783\,82469\,54\times 10^5$ in exact agreement at this level of precision.

In Table 1 we show the absolute relative error in the expansion of $S_\mu^-(a;\la)$ based on the algebraic component
\bee\label{e41}
S_\mu^-(a;\la)\sim \frac{1}{2a^{2\mu}}+H_\mu^-(a;\la)
\ee
for different values of $a$ and truncation index $k$ in the expansion of $H_\mu^-(a;\la)$ in (\ref{e24}). The final row shows the value of $S_\mu^-(a;\la)$ obtained by high-precision evaluation of (\ref{e13}). Table 2 shows the absolute relative error for complex $a$ when $a=6e^{i\phi}$ for different $\phi$, $\mu$ and $\la$.

\begin{table}[t]
\caption{\footnotesize{The absolute relative error in the computation of $S_{\mu}^-(a;\la)$ from (\ref{e41}) for different $a$ and truncation index $k$ in the asymptotic expansion $H_\mu^-(a;\la)$ when $\lambda=1$ and $\mu=\fs$.}}
\begin{center}
\begin{tabular}{|c|l|l|l|}
\hline
&&&\\[-0.3cm]
\mcol{1}{|c|}{$k$} & \mcol{1}{c|}{$a=6$} & \mcol{1}{c|}{$a=8$}& \mcol{1}{c|}{$a=10$} \\
[.1cm]\hline
&&&\\[-0.25cm]
0 & $1.775\times 10^{-3}$ & $4.859\times 10^{-4}$ & $6.275\times 10^{-4}$ \\
1 & $5.148\times 10^{-5}$ & $1.593\times 10^{-5}$ & $6.455\times 10^{-6}$ \\
2 & $2.681\times 10^{-6}$ & $4.738\times 10^{-7}$ & $1.233\times 10^{-7}$ \\
4 & $5.156\times 10^{-8}$ & $1.959\times 10^{-9}$ & $1.713\times 10^{-10}$ \\
6 & $1.278\times 10^{-8}$ & $2.411\times 10^{-10}$ & $1.000\times 10^{-11}$ \\
8 & $3.294\times 10^{-9}$ & $3.834\times 10^{-12}$ & $7.940\times 10^{-14}$ \\
[.1cm]\hline
&&&\\[-0.25cm]
$S_{\mu}^-(a;\la)$ & $1.22060\times 10^{-1}$ & $9.14725\times 10^{-2}$ & $7.31518\times 10^{-2}$\\
 [.2cm] \hline
\end{tabular}
\end{center}
\end{table}

\begin{table}[t]
\caption{\footnotesize{The absolute relative error in the computation of $S_{\mu}^-(a;\la)$ from (\ref{e41}) for different $\mu$ and $\la$ when $a=6e^{\pi i\phi}$ and truncation index $k=8$ in the asymptotic expansion $H_\mu^-(a;\la)$.}}
\begin{center}
\begin{tabular}{|c|c|c|c|}
\hline
&&&\\[-0.3cm]
\mcol{1}{|c|}{$\phi$} & \mcol{1}{c|}{$\mu=\f{1}{4},\ \la=\fs$} & \mcol{1}{c|}{$\mu=\f{3}{4},\ \la=\f{3}{2}$}& \mcol{1}{c|}{$\mu=\f{1}{3},\ \la=\f{1}{5}$} \\
[.1cm]\hline
&&&\\[-0.25cm]
0    & $4.497\times 10^{-9}$ & $4.157\times 10^{-10}$& $1.006\times 10^{-8}$ \\
0.10 & $8.383\times 10^{-9}$ & $2.293\times 10^{-9}$ & $2.798\times 10^{-8}$ \\
0.20 & $6.178\times 10^{-8}$ & $1.005\times 10^{-8}$ & $3.321\times 10^{-7}$ \\
0.30 & $2.088\times 10^{-6}$ & $1.098\times 10^{-7}$ & $1.667\times 10^{-5}$ \\
0.40 & $2.615\times 10^{-4}$ & $5.917\times 10^{-6}$ & $2.698\times 10^{-3}$ \\
[.1cm]\hline
\end{tabular}
\end{center}
\end{table}

The algebraic component of $S_\mu^+(a;\la)$ is given by
\bee\label{e42}
S_\mu^+(a;\la)\sim \frac{1}{2a^{2\mu}}+J_\mu(a;\la)+H_\mu^+(a;\la),
\ee
where the expansion of $H_\mu^+(a;\la)$ is given by (\ref{e32}). The expansion of $J_\mu(a;\la)$ can be obtained from that of ${\bf K}_\nu(z)$ given in \cite[(11.6.1)]{DLMF} to yield the asymptotic expansion
\[J_\mu(a;\la)\sim \frac{a^{1-2\mu}}{2} \sum_{k=0}^\infty \frac{(-1)^k (\fs)_k (\mu)_k}{(\fs\la a)^{2k+1}}\]
for $|a|\to\infty$ in $|\arg\,a|<\fs\pi$ (with $\la$ bounded away from zero).
The absolute relative error in the expansion of $S_\mu^+(a;\la)$ for different values of $a$ and common truncation index $k$ in the expansions of $H_\mu^+(a;\la)$ and $J_\mu(a;\la)$ are presented in Table 3, where the final row shows the value of $S_\mu^+(a;\la)$.
\begin{table}[t]
\caption{\footnotesize{The absolute relative error in the computation of $S_{\mu}^+(a;\la)$ from (\ref{e42}) for different $a$ and truncation index $k$ in the asymptotic expansions of $H_\mu^+(a;\la)$ and $J_\mu(a;\la)$ when $\lambda=1$ and $\mu=\f{1}{4}$.}}
\begin{center}
\begin{tabular}{|c|l|l|l|}
\hline
&&&\\[-0.3cm]
\mcol{1}{|c|}{$k$} & \mcol{1}{c|}{$a=10$} & \mcol{1}{c|}{$a=15$}& \mcol{1}{c|}{$a=20$} \\
[.1cm]\hline
&&&\\[-0.25cm]
0 & $2.959\times 10^{-3}$ & $1.358\times 10^{-3}$ & $7.736\times 10^{-4}$ \\
1 & $1.991\times 10^{-4}$ & $4.293\times 10^{-5}$ & $1.408\times 10^{-5}$ \\
2 & $3.864\times 10^{-5}$ & $3.962\times 10^{-6}$ & $7.525\times 10^{-7}$ \\
3 & $1.485\times 10^{-5}$ & $7.268\times 10^{-7}$ & $8.054\times 10^{-8}$ \\
4 & $9.491\times 10^{-6}$ & $2.214\times 10^{-7}$ & $1.433\times 10^{-8}$ \\
5 & $9.129\times 10^{-6}$ & $1.010\times 10^{-7}$ & $3.817\times 10^{-9}$ \\
[.1cm]\hline
&&&\\[-0.25cm]
$S_{\mu}^+(a;\la)$ & $4.98789\times 10^{-1}$ & $4.07911\times 10^{-1}$ & $3.53467\times 10^{-1}$\\
 [.2cm] \hline
\end{tabular}
\end{center}
\end{table}

\vspace{0.6cm}

\begin{center}
{\bf 4.\ Concluding remarks}
\end{center}
\setcounter{section}{4}
\setcounter{equation}{0}
\renewcommand{\theequation}{\arabic{section}.\arabic{equation}}
We have obtained a representation of the Mathieu-exponential sums $S_\mu^\pm(a;\la)$ in (\ref{e11}) in the form of an infinite sum of modified Bessel functions of complex argument together with a contribution that results from integration between the branch points $\pm ia$. For large $a$ in the sector $|\arg\,a|<\fs\pi$, this last contribution possesses an algebraic-type asymptotic expansion. When $\la=0$, our result in the case of $S_\mu^-(a;\la)$ reduces to the expression given by Olver stated in (\ref{e14}).

A problem with the representation (\ref{e26}) is that it has been derived only when $0\leq\mu<1$. It is not obvious how this result can be analytically continued into $\mu\geq 1$, as is the case in (\ref{e14}). One approach might be to exploit the fact that
\[S_{\mu+1}(a;\la)=\sum_{n=0}^\infty \frac{(-1)^n e^{-\la n}}{(n^2+a^2)^{\mu+1}}=-\frac{1}{2\mu a} \frac{\partial}{\partial a} S_\mu^-(a;\la)\qquad (0<\mu<1).\]
Differentiation of the right-hand side of (\ref{e26}) using the series expansion (\ref{e22}) for $H_\mu^-(a;\la)$ and the properties of the $K$-Bessel function shows that its value is given by (\ref{e26}) with $\mu$ replaced by $\mu+1$. This enables us to extend the representation to $1<\mu<2$, and by continuation to higher ranges of $\mu$.
The representations when $\mu=1, 2, \ldots ,5$ are displayed in (\ref{e16}). Thus, the algebraic part of the expansion is given by
\[S_\mu^-(a;\la)\sim \frac{1}{2a^{2\mu}}+H_\mu^-(a;\la)\sim\frac{1}{2a^{2\mu}}+\frac{1}{a^{2\mu}}\sum_{k=0}^\infty \frac{(\mu)_k B_k}{k! a^{2k}}\qquad (\mu\geq0)\]
as $|a|\to\infty$ in $|\arg\,a|<\fs\pi$. When $\mu=0$, we note that the right-hand side of the above expression reduces to $e^\la/(1+e^\la)$ as stated in (\ref{e15}).
Similar considerations can be brought to bear on the sum $S_\mu^+(a;\la)$.

Finally, we observe that when $\la=0$, the quantity $J_\mu(a;\la)$ defined in (\ref{e31}) reduces to
\[J_\mu(a;0)=\frac{\sqrt{\pi}\,\g(\mu-\fs)}{2a^{2\mu-1}\g(\mu)}.\]
Since $H_\mu^+(a;0)\equiv 0$, we consequently find that
\[\sum_{n=0}^\infty \frac{e^{-\la n}}{(n^2+a^2)^\mu}=\frac{1}{2a^{2\mu}}+\frac{\sqrt{\pi}\, \g(\mu-\fs)}{2a^{2\mu-1}\g(\mu)}+\frac{2^{\frac{3}{2}-\mu}\sqrt{\pi}}{a^{2\mu-1} \g(\mu)} \sum_{k=0}^\infty \frac{K_{\frac{1}{2}-\mu}((2k+2)\pi a)}{((2k+2)\pi a)^{\frac{1}{2}-\mu}}\]
for $\mu>0$ and $|\arg\,a|<\fs\pi$, which complements Olver's result in (\ref{e14}).

\end{document}